\documentclass{gtart_h}

\def\ifplaintex{\expandafter\ifx\csname documentclass\endcsname\relax}

\def\gtp{{\mathsurround=0pt\it $\cal G\mskip-2mu$eometry \&\ 
$\cal T\!\!$opology $\cal P\!$ublications}}  

\def\recd{{\small Received:\qua\receiveddate\ifx\reviseddate\relax
\else\qquad Revised:\qua\reviseddate\fi\par}} 

\def\lognumber#1{\def\thelognumber{#1}}
\def\volumenumber#1{\def\thevolumenumber{#1}}
\def\volumeyear#1{\def\thevolumeyear{#1}}
\def\papernumber#1{\def\thepapernumber{#1}}
\def\pagenumbers#1#2{\def\startpage{#1}\def\finishpage{#2}}
\def\published#1{\def\publishdate{#1}}

\def\received#1{\def\receiveddate{#1}}
\def\revised#1{\def\reviseddate{#1}}
\def\accepted#1{\def\accepteddate{#1}}

\def\asciiaddress#1{\def\theasciiaddress{#1}}
\def\asciiemail#1{\def\theasciiemail{#1}}

\let\\\par\let\thelognumber\relax\let\thevolumenumber\relax
\let\thepapernumber\relax\let\thevolumeyear\relax\let\startpage\relax
\let\finishpage\relax\let\publishdate\relax\let\receiveddate\relax
\let\reviseddate\relax\let\accepteddate\relax\let\theasciititle\relax
\let\theasciiauthors\relax\let\theasciiaddress\relax
\let\theasciiabstract\relax

\let\theasciiemail\relax

\ifplaintex
\font\logobig=cmssbx10 scaled 3836
\font\logomed=cmssbx10 scaled 2557
\else
\font\logobig=cmssbx10 scaled 4200
\font\logomed=cmssbx10 scaled 2800
\fi

\long\def\makeagttitle{   
\count0=\startpage
\agt\hfill      
\hbox to 45truept{\vbox to 0pt{\vglue -13truept{\logomed A\kern -.37em{\logobig 
T}\kern -.38em G}\vss}\hss}
\break
{\small Volume \thevolumenumber\ (\thevolumeyear)
\startpage--\finishpage\nl
Published: \publishdate}

\vglue .25truein

{\parskip=0pt\leftskip 0pt plus
1fil\def\\{\par\smallskip}{\Large\bf\thetitle}\par\medskip} \vglue
0.05truein

{\parskip=0pt\leftskip 0pt plus 1fil\def\\{\par}{\sc\theauthors}
\par\medskip}%
 
\vglue 0.03truein 

{\small\leftskip 25truept\rightskip 25truept{\bf Abstract}\stdspace\theabstract

{\bf AMS Classification}\stdspace\theprimaryclass
\ifx\thesecondaryclass\relax\else; \thesecondaryclass\fi\par
{\bf Keywords}\stdspace \thekeywords\par}\vglue 7truept

}   

\ifplaintex
\font\phead=cmsl9 scaled 950
\font\pnum=cmbx10 scaled 913
\font\pfoot=cmsl9 scaled 950
\headline{\vbox to 0pt{\vskip -4.5mm\line{\small\phead\ifnum
\count0=\startpage ISSN 1472-2739 (on-line) 1472-2747 (printed)
\hfill {\pnum\folio}\else\ifodd\count0\def\\{ }%
\ifx\theshorttitle\relax\thetitle\else\theshorttitle\fi\hfill{\pnum\folio}
\else\def\\{ and }{\pnum\folio}\hfill\ifx\theshortauthors\relax\theauthors
\else\theshortauthors\fi\fi\fi}\vss}}
\footline{\vbox to 0pt{\vglue 0mm\line{\small\pfoot\ifnum\count0=\startpage
\copyright\ \gtp\hfill\else
\agt, Volume \thevolumenumber\ (\thevolumeyear)\hfill\fi}\vss}}
\else
\font=cmsl9 scaled 1050
\font\lnum=cmbx10 
\font=cmsl9 scaled 1050
\makeatletter
\def\@oddhead{{\small\lhead\ifnum\count0=\startpage ISSN 1472-2739 
(on-line) 1472-2747 (printed)\hfill {\lnum\number\count0}\else\ifodd\count0
\def\\{ }\ifx\theshorttitle\relax \thetitle \else\theshorttitle\fi\hfill
{\lnum\number\count0}\else\def\\{ and }{\lnum\number\count0}
\hfill\ifx\theshortauthors\relax 
\theauthors\else\theshortauthors\fi\fi\fi}}\def\@evenhead{\@oddhead}
\def\@oddfoot{\small\lfoot\ifnum\count0=\startpage\copyright\ \gtp\hfill\else
\agt, Volume \thevolumenumber\ (\thevolumeyear)\hfill\fi}
\def\@evenfoot{\@oddfoot}
\makeatother
\fi
\let\maketitlepage\makeagttitle
\let\makeshorttitle\maketitlepage
\let\maketitle\maketitlepage

\newwrite\gtoutfile
\long\gdef\makeheadfile{  
{\def\\{, }\def\s{ }
\immediate\openout\gtoutfile head.xxx
\immediate\write\gtoutfile{Proxy-for: \ifx\theasciiauthors\relax
\theauthors\else\theasciiauthors\fi\s<\ifx\theasciiemail\relax\theemail\else\theasciiemail\fi>}
\immediate\write\gtoutfile{\noexpand\\}
\immediate\write\gtoutfile{Authors: \ifx\theasciiauthors\relax
\theauthors\else\theasciiauthors\fi}
{\def\\{ }\immediate\write\gtoutfile{Title: \ifx\theasciititle\relax
\thetitle\else\theasciititle\fi}}
\immediate\write\gtoutfile{Subj-class: GT or SG, GR etc}
\immediate\write\gtoutfile{MSC-class: \theprimaryclass\ifx\thesecondaryclass\relax\else, \thesecondaryclass\fi}
\immediate\write\gtoutfile{Journal-ref: Algebr. Geom. Topol. \thevolumenumber\s
(\thevolumeyear) \startpage-\finishpage}
\immediate\write\gtoutfile{Comments: Published by Algebraic and
Geometric Topology at}
\immediate\write\gtoutfile{\s\s\s  http://www.maths.warwick.ac.uk/agt/AGTVol\thevolumenumber/agt-\thevolumenumber-\thepapernumber.abs.html}
\immediate\write\gtoutfile{\noexpand\\}
\immediate\write\gtoutfile{}
\ifx\theasciiabstract\relax
\immediate\write\gtoutfile{\theabstract}\else
\immediate\write\gtoutfile{\theasciiabstract}\fi
\immediate\write\gtoutfile{}
\immediate\write\gtoutfile{\noexpand\\}
\immediate\write\gtoutfile{}
\immediate\closeout\gtoutfile}}  

\def\maketitlepage{\makeagttitle\makeheadfile}
\let\makeshorttitle\maketitlepage
\let\maketitle\maketitlepage

\lognumber{6}
\volumenumber{5}
\volumeyear{2005}
\papernumber{6}
\pagenumbers{107}{118}
\received{2 February 2004} 
\revised{3 February 2005}
\accepted{6 June 2004}
\published{7 February 2005}

\usepackage{amssymb,amsmath,amscd,epsfig}

\newtheorem{thm}{Theorem}[section]    
\newtheorem{lem}[thm]{Lemma}          
\newtheorem{proposition}[thm]{Proposition}
\theoremstyle{definition}

\newcommand{\slbqone}{\raisebox{-26pt}{\mbox{}\hspace{1pt}
                  \epsfig{file=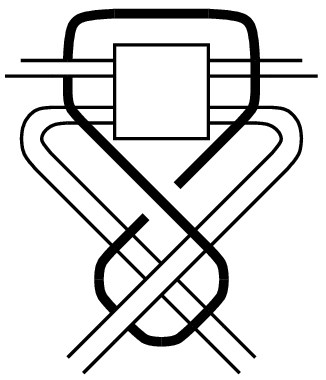,height=54pt}
                  \hspace{1pt}\mbox{}}}

\newcommand{\bk}{\mbox{\raisebox{-16pt}{\mbox{}\hspace{1pt}
                 \epsfig{file=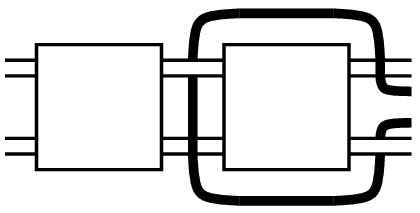,height=37pt}}
                 \hspace{-70pt}{\tiny{$q-k$}}\hspace{17pt}
		 {\tiny{$k$}}\hspace{22pt}\mbox{}}}

\newcommand{\yminus}{\mbox{\raisebox{-16pt}{\mbox{}\hspace{1pt}
                 \epsfig{file=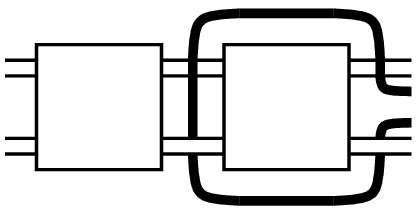,height=37pt}}
                 \hspace{-63pt}{\tiny{$q$}}\hspace{26pt}
		 {\tiny{$0$}}\hspace{22pt}\mbox{}}}

\newcommand{\bmo}{\mbox{\raisebox{-16pt}{\mbox{}\hspace{1pt}
                  \epsfig{file=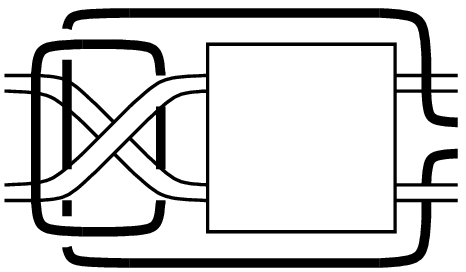,height=36pt}}
		  \hspace{-27pt}\smst{k-1}\hspace{28pt}\mbox{}}}
\newcommand{\bmi}{\mbox{\raisebox{-16pt}{\mbox{}\hspace{1pt}
                  \epsfig{file=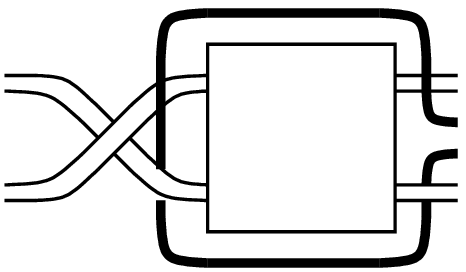,height=36pt}}
		  \hspace{-27pt}\smst{k-1}\hspace{28pt}\mbox{}}}
\newcommand{\bmii}{\mbox{\raisebox{-16pt}{\mbox{}\hspace{1pt}
                  \epsfig{file=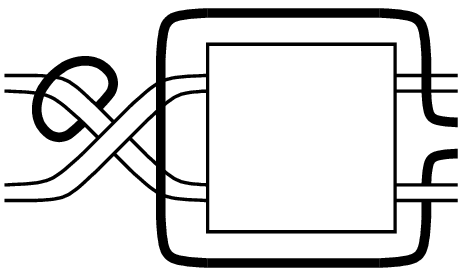,height=36pt}}
		  \hspace{-27pt}\smst{k-1}\hspace{28pt}\mbox{}}}
\newcommand{\bmiii}{\mbox{\raisebox{-16pt}{\mbox{}\hspace{1pt}
                  \epsfig{file=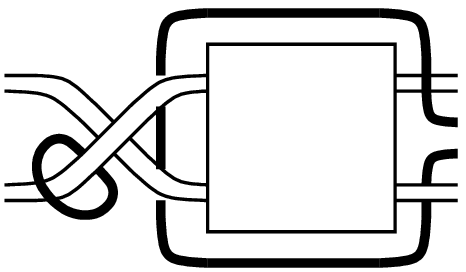,height=36pt}}
		  \hspace{-27pt}\smst{k-1}\hspace{28pt}\mbox{}}}
\newcommand{\bmiv}{\mbox{\raisebox{-16pt}{\mbox{}\hspace{1pt}
                  \epsfig{file=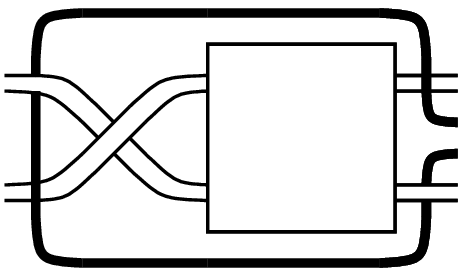,height=36pt}}
		  \hspace{-26pt}\smst{k-1}\hspace{28pt}\mbox{}}}

\newcommand{\ybo}{\mbox{\raisebox{-23pt}{\mbox{}\hspace{1pt}
                  \epsfig{file=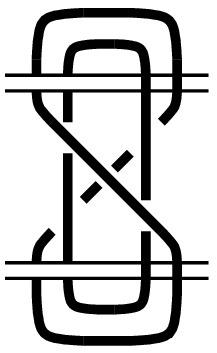,height=50pt}}
		  \hspace{-8pt}\raisebox{-2pt}{\smst{\mbox{}^m}}\mbox{}}}

\newcommand{\sbi}{\raisebox{-26pt}{\mbox{}\hspace{1pt}
                  \epsfig{file=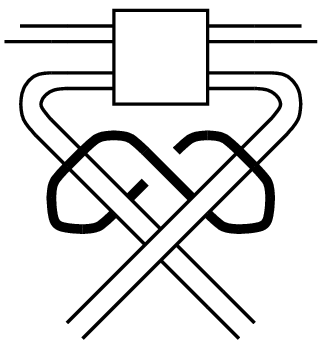,width=50pt}
                  \hspace{1pt}\mbox{}}}
\newcommand{\sbii}{\raisebox{-26pt}{\mbox{}\hspace{1pt}
                  \epsfig{file=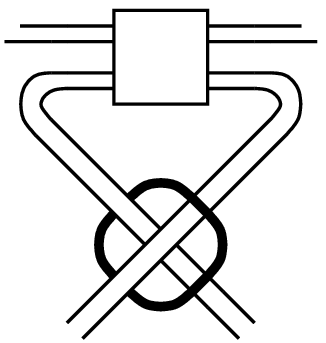,width=50pt}
                  \hspace{1pt}\mbox{}}}
\newcommand{\xto}{\raisebox{-26pt}{\mbox{}\hspace{1pt}
                  \epsfig{file=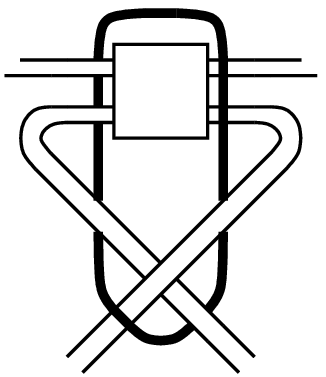,width=50pt}
                  \hspace{1pt}\mbox{}}}
\newcommand{\xti}{\raisebox{-26pt}{\mbox{}\hspace{1pt}
                  \epsfig{file=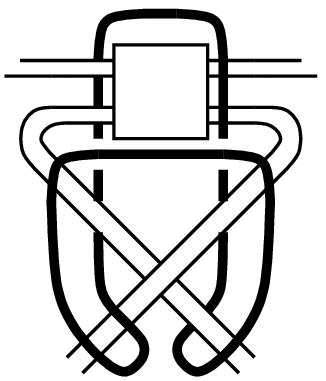,width=50pt}
                  \hspace{1pt}\mbox{}}}
\newcommand{\xtii}{\raisebox{-26pt}{\mbox{}\hspace{1pt}
                  \epsfig{file=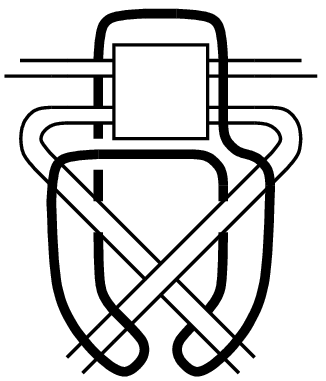,width=50pt}
                  \hspace{1pt}\mbox{}}}
\newcommand{\xtiii}{\raisebox{-26pt}{\mbox{}\hspace{1pt}
                  \epsfig{file=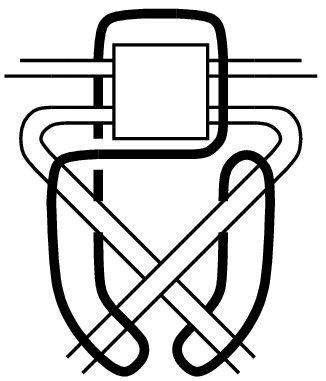,width=50pt}
                  \hspace{1pt}\mbox{}}}
\newcommand{\xtiv}{\raisebox{-26pt}{\mbox{}\hspace{1pt}
                  \epsfig{file=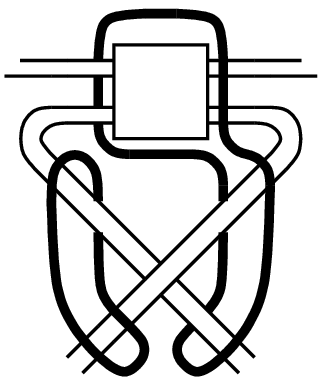,width=50pt}
                  \hspace{1pt}\mbox{}}}
\newcommand{\xtv}{\raisebox{-26pt}{\mbox{}\hspace{1pt}
                  \epsfig{file=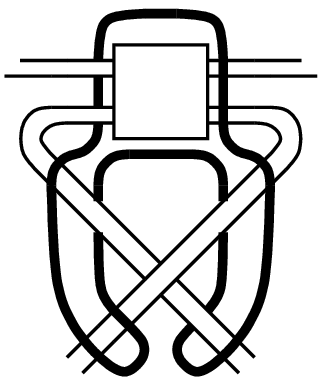,width=50pt}
                  \hspace{1pt}\mbox{}}}
\newcommand{\xtvi}{\raisebox{-26pt}{\mbox{}\hspace{1pt}
                  \epsfig{file=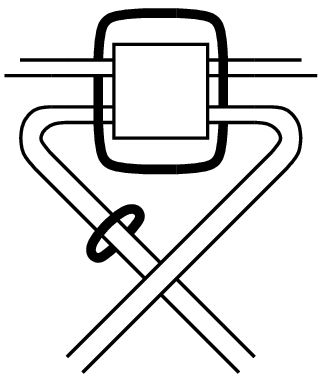,width=50pt}
                  \hspace{1pt}\mbox{}}}
\newcommand{\xtvii}{\raisebox{-26pt}{\mbox{}\hspace{1pt}
                  \epsfig{file=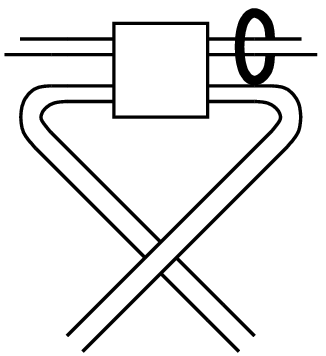,width=50pt}
                  \hspace{1pt}\mbox{}}}
\newcommand{\xtviii}{\raisebox{-26pt}{\mbox{}\hspace{1pt}
                  \epsfig{file=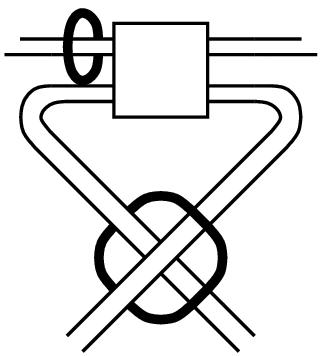,width=50pt}
                  \hspace{1pt}\mbox{}}}

\newcommand{\ukx}{\raisebox{-11pt}{\mbox{}\hspace{1pt}
                  \epsfig{file=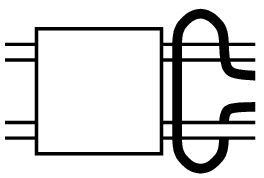,height=27pt}
                  \hspace{1pt}\mbox{}}}

\newcommand{\ukz}{\raisebox{-11pt}{\mbox{}\hspace{1pt}
                  \epsfig{file=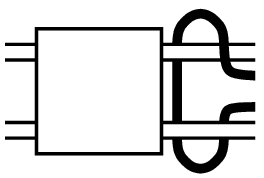,height=27pt}
                  \hspace{1pt}\mbox{}}}
\newcommand{\uku}{\raisebox{-11pt}{\mbox{}\hspace{1pt}
                  \epsfig{file=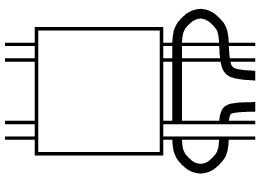,height=27pt}
                  \hspace{1pt}\mbox{}}}

\newcommand{\uti}{\raisebox{-16pt}{\mbox{}\hspace{1pt}
                  \epsfig{file=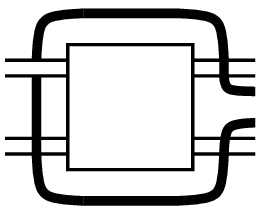,height=37pt}
                  \hspace{1pt}\mbox{}}}
\newcommand{\utii}{\raisebox{-13pt}{\mbox{}\hspace{1pt}
                  \epsfig{file=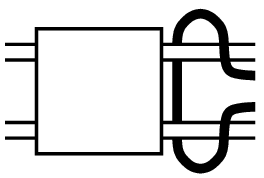,height=32pt}
                  \hspace{1pt}\mbox{}}}
\newcommand{\utiii}{\raisebox{-13pt}{\mbox{}\hspace{1pt}
                  \epsfig{file=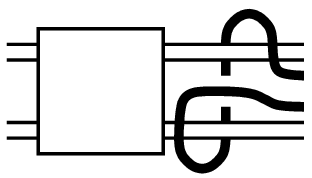,height=32pt}
                  \hspace{1pt}\mbox{}}}

\newcommand{\lcr}{\raisebox{-5pt}{\mbox{}\hspace{1pt}
                  \epsfig{file=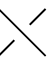}\hspace{1pt}\mbox{}}}
\newcommand{\ift}{\raisebox{-5pt}{\mbox{}\hspace{1pt}
                  \epsfig{file=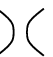}\hspace{1pt}\mbox{}}}
\newcommand{\zer}{\raisebox{-5pt}{\mbox{}\hspace{1pt}
                  \epsfig{file=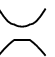}\hspace{1pt}\mbox{}}}

\newcommand{\smst}[1]{\makebox[0pt]{\scriptsize{$#1$}}}
\newcommand{\lmk}{\noindent\mbox{}\hfill}
\newcommand{\rmk}{\hfill\mbox{}\par}

\newcommand{\cL}{\mathcal L}
\newcommand{\cB}{\mathcal B}

\newcommand{\cR}{\mathcal R}

\newcommand{\bZ}{\mathbb Z}

\newcommand{\bC}{\mathbb C}

\newcommand{\g}{SL_2(\bC)}

\begin{document}
\title{The Kauffman bracket skein module\\of a twist knot exterior}

\authors{Doug Bullock\\Walter Lo Faro}
\address{Department of Math, Boise State University\\Boise, ID 83725, USA}

\secondaddress{Department of Math, University of Wisconsin\\Stevens
Point, WI 54481, USA}
\asciiaddress{Department of Math, Boise State University\\Boise, ID 83725, 
USA\\and\\Department of Math, University of Wisconsin\\Stevens
Point, WI 54481, USA}

\gtemail{\mailto{bullock@math.boisestate.edu}{\rm\qua 
and\qua}\mailto{Walter.LoFaro@uwsp.edu}}
\asciiemail{bullock@math.boisestate.edu, Walter.LoFaro@uwsp.edu}

\begin{abstract} 
We compute the Kauffman bracket skein module of the complement of a
twist knot, finding that it is free and infinite dimensional.  The
basis consists of cables of a two-component link, one component of
which is a meridian of the knot.  The cabling of the meridian can be
arbitrarily large while the cabling of the other component is limited
to the number of twists.
\end{abstract}

\primaryclass{57M27}                
\secondaryclass{57M99}              
\keywords{Knot, link, skein module, Kauffman bracket}                    
\maketitle

\section{Introduction}

At first glance, and in original intent \cite{P1}, the Kauffman
bracket skein module is a formal extension of the Kauffman bracket
polynomial to an arbitrary 3-manifold.  As Kauffman's polynomial (for
framed links in $S^3$) is equivalent to the Jones polynomial (for
oriented links in $S^3$), one may think of the skein module as a
generalization of the Jones polynomial.  More recently the module has
taken on a different significance: it is now seen as a deformation of
the $\g$-characters of the fundamental group
\cite{isomorphism,BFK,PS}.  Using this interpretation of the skein
module of a knot exterior, Frohman, Gelca and the second author here
constructed a quantum version the $A$-polynomial \cite{FGL}.  This is
related back to the Jones polynomial \cite{T,Y} (not simply by
generality) and has implications for the hyperbolic volume conjecture
\cite{K,MM}.  Despite all this, there have as yet been no computations
of Kauffman bracket skein modules for hyperbolic manifolds.

Early computations \cite{P1} depended on an $I$-bundle structure for
the manifold, since projection along the $I$ factor gave a natural
mechanism for controlling complexity.  The only other successful
method \cite{skein,HP1,HP2} has been to consider the effect of adding
a single 2-handle to a handlebody.  This creates a presentation with
fairly simple generators (any basis for the module of the handlebody),
but having an unwieldy set of relations.  Eliminating redundant
relations is the most difficult part of the task.  What is needed is
an effective method of creating relations among relations, or {\em
syzygies}, and then keeping track of which relations can be removed.

This has been managed for all genus one manifolds \cite{HP1,HP2}, and
for $(2,q)$-torus knot exteriors \cite{skein}.  In principle, the
combinatorics ought to be accessible for genus two manifolds with
toral boundary (one added handle), but the computations are quite
daunting in practice.  Even for $(2,q)$-torus knots, the trick was
managed only with help from a particularly nice basis.

The innovation in this paper is a simpler method of keeping track of
the relations.  We use the established handle addition technique, but
we twist the handlebody instead of the handle, which simplifies many
bracket computations.  Our viewpoint also leads to a comfortable and
practical method for producing syzygies that reduce the initial
presentation to a simple basis.

\section{The theorem}

Let $M$ be an orientable 3-manifold.  A framed link in $M$ is an
embedding of a disjoint collection of annuli into $M$.  Framed links
are depicted by link diagrams showing the cores of an annuli lying
flat in the projection plane (i.e.\ with blackboard framing).

Two framed links in $M$ are equivalent if there is an isotopy of $M$
taking one to the other.  Let $\cL_M$ denote the set of equivalence
classes of framed links in $M$, including the empty link.  With $R =
\bZ [t^{\pm 1}]$, form a free module $R \cL_M$ with basis $\cL_M$.
Define $S(M)$ to be the smallest submodule of $R \cL_M$ containing all
expressions of the form $\displaystyle{\lcr - t \zer - t^{-1} \ift}$
and $\bigcirc + t^2 + t^{-2}$, where the framed links in each
expression are identical outside balls pictured in the diagrams.  The
Kauffman bracket skein module $K(M)$ is the quotient
\[ R \cL_M / S(M). \]

\begin{figure}
\lmk\epsfig{file=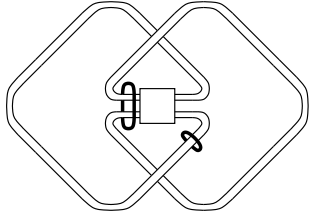,height=81pt}\rmk
\caption{Exterior of a $q$-twist knot and the link $xy$}
\label{knot}
\end{figure}

\begin{figure}
\lmk\epsfig{file=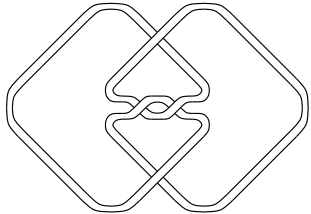,height=81pt}\rmk
\caption{Figure-8 knot as a 2-twist knot}
\label{example}
\end{figure}

A $q$-twist knot (right-handed, if not amphichiral) is the alternating
knot formed by inserting a left-to-right string of $q$ half-twists
into the coupon in Figure \ref{knot}.  The 2-twist knot in Figure
\ref{example}, for example, is the familiar figure-8 knot.  Let $M_q$
be the twist knot exterior, and denote by $xy$ the 0-framed, two
component link also pictured in Figure \ref{knot}.  The meridian is
$x$ and the other component is $y$.  In general, $x^ly^m$ denotes the
cable of $xy$ consisting of $l$ parallel copies of $x$ and $m$
parallel copies of $y$.  The exponents run over non-negative integers
and 1 denotes the empty link.

\begin{thm}
$K(M_q)$ is free with basis $\{x^ly^m \;|\; m \leq q\}$
\end{thm}

\section{Initial presentation}

\begin{figure}
\lmk\epsfig{file=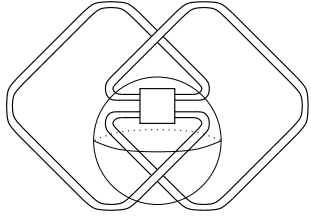,height=81pt}\rmk
\caption{Decomposing sphere $S$ in $M_q$}
\label{decomposed}
\end{figure}

\begin{figure}
\lmk\epsfig{file=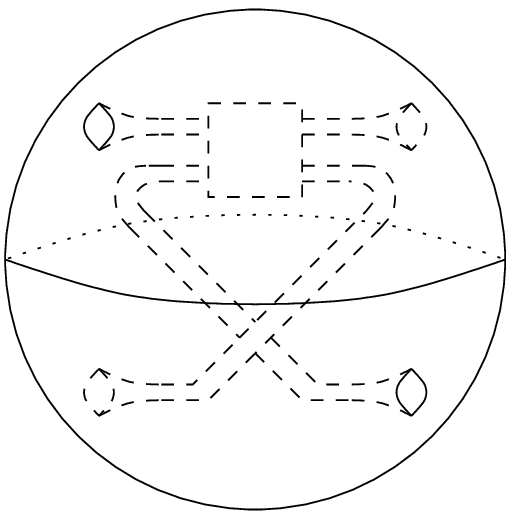,height=72pt}
\hfill\epsfig{file=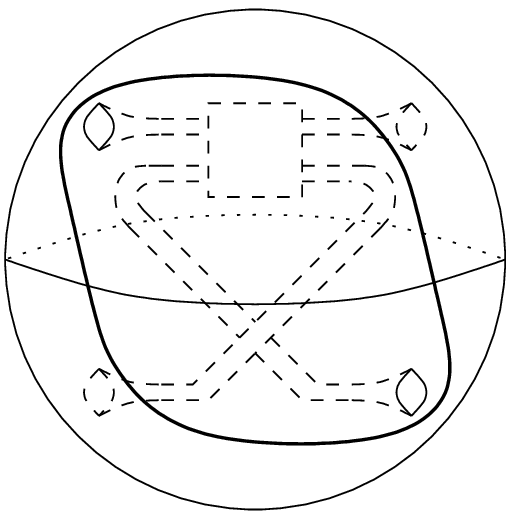,height=72pt}
\hfill\epsfig{file=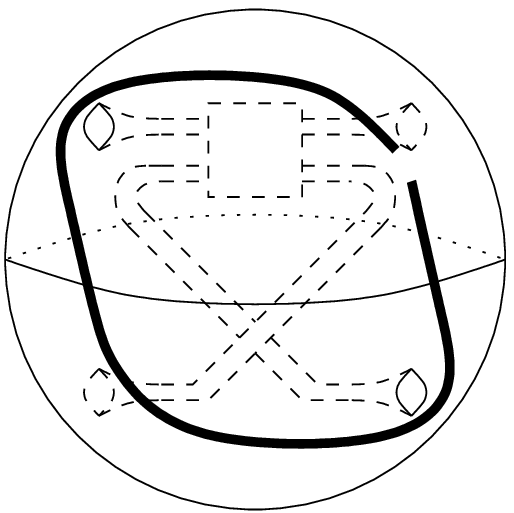,height=72pt}\rmk
\lmk\makebox[72pt]{(a)}
\hfill\makebox[72pt]{(b)}
\hfill\makebox[72pt]{(c)}\rmk
\caption{(a) $H$\qua  (b) Core of the attaching annulus\qua (c) The core as
a relative link}
\label{newcore}
\end{figure}

The knot exterior $M_q$ decomposes into a pair of genus two
handlebodies glued along the 4-punctured sphere $S$ shown in Figure
\ref{decomposed}.  Let $H$ be the closure of the component of $M_q-S$
containing the coupon.  Figure \ref{newcore}(a) depicts $H$, slightly
deformed so that the upper left and lower right punctures are in the
foreground.

The portions of the knot outside $H$ are parallel to a pair of arcs in
$S$ that cut it into an annulus.  Therefore, $M_q$ is homeomorphic to
$H$ with a 2-handle attached along this annulus.  Its core is shown in
Figure \ref{newcore}(b).

There is a standard argument \cite{skein,surgery,HP1,HP2,W} that says
$K(M_q)$ is $K(H)$ modulo skeins differing by slides across the
2-handle.  We find this language to be a little imprecise, so we will
rephrase it in terms of relative skeins.  Suppose that the core of the
attaching annulus is given the blackboard framing in $S$.  We cut out
a very small bit of this curve, as indicated in Figure
\ref{newcore}(c), leaving a framed arc whose ends are a pair of framed
points in $\partial H$.  Following \cite{Lo,P1}, let $K_1(H)$ be the
skein module of $H$ relative to those two framed points.

Let $L$ be a relative link in $H$.  Since the ends of $L$ are very
close together we can unambiguously define the {\em completion} of $L$
to be the result of gluing its ends together.  The {\em slide} of $L$
is formed by gluing its ends to the cut open core of the attaching
annulus.  Completion and slide are denoted by $c(L)$ and $s(L)$.  Let
$r(L) = c(L) - s(L)$ and extend linearly to $r \co K_1(H) \to K(H)$.
The image of $r$ is the set of all possible relations in $K(H)$
induced by handle slides.  Therefore,
\[K(M_q) = K(H) / r(K_1(H)).\]

Since $K(H)$ is free, the quotient provides a presentation of
$K(M_q)$.  Any basis for $K(H)$ serves as a generating set.  For
relations, choose generators for $K_1(H)$, apply $r$, and express
everything in terms of the basis of $K(H)$.  The more efficient your
generating set for $K_1(H)$, the more efficient your presentation.
However, even a basis for $K_1(H)$ yields unnecessary relations.
Computing $K(M_q)$ thus becomes a search for all relations among the
relations in this presentation.

We need to fix a basis for $K(H)$.  Let $x$ and $y$ be the knots in
Figure \ref {knot}, but only up to isotopy in $H$.  Let $z$ be a
meridian that is not isotopic to $x$ in $H$.  As usual, $x$, $y$ and
$z$ are 0-framed.  The set of cables, $\cB = \{x^ly^mz^n\}$, is a
basis for $K(H)$ \cite{P1}.

We also need to fix generators of $K_1(H)$, but first some notes on
multiplicative notation for (possibly relative) links in $H$.

\begin{itemize}
\item 
The notation is commutative and associative.
\item 
$x^l$, $y^m$ and $z^n$ denote cables.
\item 
If $L$ is a (possibly relative) link then $x^lL$ means the union of
$L$, pushed away from the knot boundary, with a cable of $x$ very near
the knot boundary.
\item 
Similarly for $z^nL$.
\item 
If $\sigma$ is a (possibly relative) skein then $x^lz^n\sigma$ is
defined by distributing $x^lz^n$ across any linear combination of
links representing $\sigma$.  This is well defined because
representatives of $\sigma$ differ by skein relations that take place
away from the knot boundary.
\item 
If $\sigma$ is written in terms of the basis $\{x^iy^jz^k\}$ then
$x^lz^n\sigma$ is just polynomial multiplication.
\item 
In general, $y^mL$ is not well defined, but there are some
specific embeddings of $L$ for which we want $y^mL$ to make sense.
These are explained below.
\end{itemize}

If $L$ is one of the relative knots
\[X = \ukx, \quad  Z = \ukz \quad \text{or} \quad U = \uku\]
then $y^mL$ denotes a copy of $y^m$ inserted into the twist coupon. If
$0 \leq k \leq q$, let
\[Y_k = \bk\] 
where the coupons contain $q-k$ and $k$ twists.  By $y^mY_k$ we mean a
cable of $y$ inserted into the coupon containing $k$ twists, even if
$k = 0$.  Lastly,
\[y^mY_{-1} = \yminus \]
with $y^m$ inserted in the coupon that contains no twists.  

\begin{lem}\label{basis}
$K_1(H)$ is generated by $\{x^ly^mz^nL\;|\; \text{$L = X$, $Y_0$, $Z$,
or $U$}\}$.\footnote{It's actually a basis, but the proof is annoying
and the result is unnecessary.}
\end{lem}

\begin{figure}
\lmk\epsfig{file=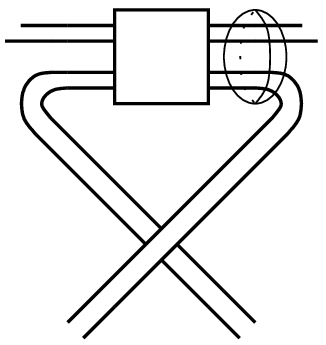,height=54pt}\rmk
\nocolon\caption{}
\label{ovalnbhd}
\end{figure}

\begin{proof}
Given any relative link in $H$, it can be isotoped into the oval
neighborhood shown in Figure \ref{ovalnbhd}.  Once there, grab the top
and bottom of the oval and twist in opposite directions a quarter turn
each.  This should make the tubes perpendicular to the page so that
\[X = \raisebox{-10pt}{\epsfig{file=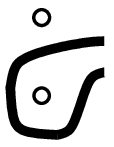,height=20pt}}\;,\quad
  Y_0 = \raisebox{-10pt}{\epsfig{file=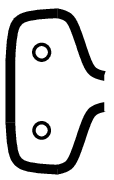,height=25pt}}\;,\quad
  Z = \raisebox{-7pt}{\epsfig{file=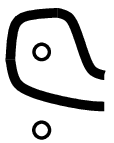,height=20pt}} \quad
  \text{and} \quad U =
  \raisebox{-7pt}{\epsfig{file=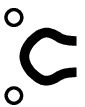,height=15pt}}\;.\] 
If not, twist the opposite way and it will.

Now resolve according to a relative version of the argument in
\cite[Lemmas 1--3]{states}.  Each term of the resolution will be a
cable of
\[\epsfig{file=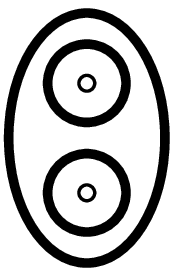,height=38pt}\]
together with one of $X$, $Y_0$, $Z$ or $U$.  The modification
introduced in \cite[Theorem 6.2]{PS} lets us force an $X$ to end up
above the cabled link and a $Z$ to end up below it.  Neither $Y_0$ nor
$U$ can become entangled.  Now we untwist the oval neighborhood,
returning $X$, $Y_0$, $Z$ and $U$ to their initial embeddings.  This
will twist the cabled link , but it can be further resolved into a
polynomial in $x$, $y$ and $z$.
\end{proof}

\section{Sufficient relations}\label{relations}

In this section we locate in $r(K_1(H))$ sufficient relations to
eliminate all but $\{x^ly^m \;|\; m \leq q\}$ from $\cB$.  It turns
out that powers of $z$ are easy to eliminate and that powers of $y$
index the complexity of other computations.  For this reason, we
introduce the notation $\sigma \sim y^m$, meaning $\sigma = t^{\pm
p}y^m$ modulo the span of $\{x^iy^jz^k\;|\;j < m\}$.

Where $x$ and $z$ are concerned, the relation submodule behaves like
an ideal.

\begin{lem}\label{ideal}
If $\sigma \in r(K_1(H))$, then $x^iz^k \sigma \in r(K_1(H))$.
\end{lem}

\begin{proof}
Suppose $L$ is a relative link in $H$.  Since $x$ and $z$ are nowhere
near the attaching annulus on $\partial H$, it's easy to see that
$x^iz^k r(L) = r(x^iz^k L)$.  This extends to all of
$K_1(H)$.
\end{proof}

In practice, you compute a relation by grinding some $r(\sigma)$ down
to a polynomial in $x$, $y$ and $z$.  Lemma \ref{ideal} then says that
any formal multiple of that relation by $x^iz^k$ is another relation.
For example, from the relation $r(y^m Z)$ we obtain a class of relations:
\begin{align}\label{xzrels}
x^{l-1}z^n r(y^mZ) &= x^{l-1}z^n (y^mz- y^mx) \\
                         &=  x^{l-1}y^mz^{n+1} - x^ly^mz^n \notag
\end{align}

Relations (\ref{xzrels}) can be used to eliminate $\{x^ly^mz^n\;|\;n >
0\}$ from the presentation of $K(M_q)$.  Powers of $y$ are more
troublesome. To eliminate $\{x^ly^m\;|\; y > q\}$, we need some new
relations.

\begin{lem}\label{bcomplete}
$c(y^m Y_k) \sim y^{m+k+1}$
\end{lem}

\begin{proof}
Induct on $k$.  If $k$ = 0, we have
\[c(y^mY_0) = -t^{-3}\ybo\]
By counting wrapping numbers in each term of the resolution one can
see that $c(y^m Y_0) \sim y^{m+1}$.  

Another wrapping number argument shows that no term of $c(y^m Y_{-1})$
has a power of $y$ larger than $m$.
	
For $k \geq 1$ consider the relation
\begin{align}\label{fundamental}
\bmo & = t^2 \bmi + t^{-2} \bmiv \\ 
     & + \bmiii + \bmii \notag
\end{align}
Insert $y^m$ into the coupon and take the closure of every term to get
\[c(y^{m+1} Y_{k-1}) = t^2 c(y^m Y_{k-2}) + t^{-2} c(y^m Y_k) + y^m
(\text{meridians}),\]
which can be solved for $y^mY_k$.
\end{proof}

\begin{lem}\label{bslide}
$s(y^m Y_q) \sim y^{q+m}$.
\end{lem}

\begin{proof}
Note first that
\begin{align}\label{slideyq}
s(Y_q) &= \sbii  = -t^3 \sbi \\
       &= -t^4c(Y_{q-1}) - t^2(\text{meridians}) \notag
\end{align}
Then insert $y^m$ into the coupon and apply Lemma \ref{bcomplete}.
\end{proof}

Lemmas \ref{bcomplete} and \ref{bslide} imply $r(y^mY_q) \sim
y^{q+m+1}$.  Extended to include powers of $x$, there relations serve
to eliminate any remaining terms of $\cB$ with $y$-degree greater than
$q$.  Hence, with
\[\cR = \{x^iz^k r(y^mZ)\} \cup \{x^i r(y^mY_q)\}\]

\begin{proposition}
The presentation $\cB$ modulo $\cR$ reduces to the free presentation
of Theorem 1.
\end{proposition}

To finish the proof of Theorem 1 we must find relations among the
relations $r(K_1(H))$ sufficient to write them all in terms of
$\cR$.  Such a relation among relations is called a {\em syzygy}.

\section{Syzygies}

Here we show that $\cR$ generates $r(K_1(H))$.  Recall that $\cR$
contains relations of the form
\begin{align*}
&r(x^ly^mz^nZ),\; \text{and}\\
&r(x^ly^mY_q)
\end{align*}
We need to show that the span of $\cR$, denoted $\langle\cR\rangle$,
contains all relations of the form 
\begin{align*}
&r(x^ly^mz^nX),\\
&r(x^ly^mz^nY_0),\;\text{and}\\
&r(x^ly^mz^nU)
\end{align*}

\begin{lem}\label{xforz}
If $L$ is any link in $H$ (or any skein in $K(H)$), then $xL-zL \in
\langle \cR \rangle$.
\end{lem}

\begin{proof}
Express $L$ in terms of the basis $\cB$ and then apply (\ref{xzrels}). 
\end{proof}

\begin{lem}\label{x_syz}
$r(x^ly^mz^nX) \in \langle \cR \rangle$
\end{lem}

\begin{proof}
Slide $X$ and resolve as
\begin{align*}
s(X) &= \xto = \xti  \\
     &= t \xtiii + t^{-1} \xtii \\ 
     &= t \xtiii + \xtv + t^{-2} \xtiv \\
     &= - t^{-2} \xtviii + \xtvii + t^{-2} \xtvi
\end{align*}
Modulo terms of the form $xL - zL$, this is
\[ s(X) = -t^{-2}xs(Y_q) + c(X) +t^{-2}xc(Y_q)\]
which is the syzygy $r(X) = -t^{-2} r(x Y_q)$.  Inserting $y^m$ into the
coupon and multiplying by $x^lz^n$ gives the syzygy 
\[r(x^ly^mz^nX) = -t^{-2} r(x^{l+1}y^mz^n Y_q)\]  
Finally, on the right hand side, convert $z$'s to $x$'s by repeated
applications of
\begin{align*}
r(x^iy^mz^kY_q) 
&= r(x^iy^mz^kY_q) - r(x^{i+1}y^mz^{k-1}Y_q) + r(x^{i+1}y^mz^{k-1}Y_q) \\
&= zr(x^iy^mz^{k-1}Y_q) - xr(x^iy^mz^{k-1}Y_q) + r(x^{i+1}y^mz^{k-1}Y_q)
\end{align*} 
This will express $r(x^ly^mz^nY_q)$ as $r(x^{l+n}y^mY_q)$ plus terms
of the form $xL - zL$.  Lemma \ref{xforz} insures that all terms are
in $\langle \cR \rangle$.
\end{proof}

\begin{lem}\label{lastsyz}
Modulo relations of the form $xL - zL$, we have the syzygy
\[r(y^mY_q) = t^4 r(y^mY_{q-1})\]
\end{lem}

\begin{proof}
Leaving obvious isotopies to the reader, 
\[t^4s(Y_{q-1}) = -t \slbqone  = - t^2 x^2 - c(Y_q) \]
Subtract this equation from Equation (\ref{slideyq}) and insert $y^m$
as usual.
\end{proof}

\begin{lem}\label{bigsyz}
For $0 \leq k \leq q$, $r(x^ly^mY_k) \in \langle \cR \rangle$.
\end{lem}

\begin{proof}
Induct downward on $k$.  If $k=q$ we are looking at $r(x^ly^mY_q)$.
If $k = q-1$, apply the syzygy from Lemma \ref{lastsyz} multiplied by
$x^l$.  

If $k \leq q-2$, apply $r$ to Equation (\ref{fundamental}) $k+1$
twists in the coupon.  Modulo terms of the from $xL - zL$, this
becomes the syzygy
\[r(yY_{k+1}) = t^2r(Y_k) + t^{-2}r(Y_{k+2}) + r(xX) + r(xZ)\] 
Solve for $r(Y_k)$, multiply by $x^l$, and insert $y^m$ in the usual
place.
\end{proof}

\begin{lem}
$r(x^ly^mz^nY_0) \in \langle \cR \rangle$.
\end{lem}

\begin{proof}
Convert $r(x^ly^mz^lY_0)$ to $r(x^{l+n}y^mY_0)$ as in the proof of
Lemma \ref{x_syz}.  Then apply Lemmas \ref{xforz} and \ref{bigsyz}.
\end{proof}

\begin{figure}
\lmk\epsfig{file=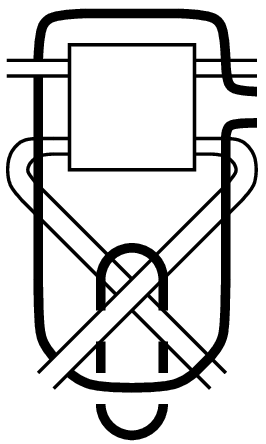,height=60pt}\rmk
\nocolon\caption{}
\label{utail}
\end{figure}

\begin{lem}
$r(x^ly^mz^nU) \in \langle \cR \rangle$.
\end{lem}

\begin{proof}
Consider the link in Figure \ref{utail}.  The relative component is
$Y_q$, and the closed component isotops into the coupon where it
resolves into some polynomial $p(x,y,z)$.  On the other hand, by
resolving the crossings as shown we obtain (modulo terms of the form
$xL-zL$)
\begin{align*}
p Y_q &= t^2 U + x \uti + x X + t^{-2} \utii \\
      &= t^2 U + x \uti + x X - t \utiii \\
      &= t^2 U + x \uti + x X - t^2 Y_0 - x X
\end{align*} 
Now apply $r$ to this equation to obtain the syzygy
\[r(p Y_q) = t^2 r(U) \pm r(xZ) -t^2 r(Y_0) \]
(The sign of $r(xZ)$ depends on the number of twists in the coupon.)
Insert $y^m$, multiply by $x^lz^n$, and solve for $r(U)$.  Except for
the term $r(x^ly^mz^nU)$, convert all $z$'s to $x$'s as usual.  The
resulting linear combination will lie in $\langle \cR \rangle$
\end{proof}

We have shown that $\langle \cR \rangle = r(K_1(H))$, so $K(M_q)$ must
be presented as in Proposition 1.

%
%
%
\Addresses\recd
%
%
%
%
\end{document}